\documentclass[10pt]{amsart}

\usepackage{verbatim}
\usepackage{eucal,url,amssymb,stmaryrd,
enumerate,amscd,
}

\usepackage{amsfonts}
\usepackage{amsmath,amsthm,amssymb,amscd,enumerate,eucal,url,stmaryrd}

\numberwithin{equation}{section}

\newtheorem{thrm}{Theorem}[section]

\newtheorem{lemma}[thrm]{Lemma}

\newtheorem{prop}[thrm]{Proposition}

\newtheorem{rmrk}[thrm]{Remark}

\def\frg{{\frak g}}
\def\frh{{\frak h}}

\begin{document}

\begin{abstract}
We construct new explicit compact supersymmetric valid solutions
with non-zero field strength, non-flat instanton and constant
dilaton to the heterotic equations of motion in dimension six. We
present balanced Hermitian structures on compact nilmanifolds in
dimension six satisfying the heterotic supersymmetry equations with
non-zero flux, non-flat instanton  and constant dilaton which obey
the three-form Bianchi identity with curvature term taken with
respect to either the Levi-Civita, the (+)-connection or the Chern
connection. Among them, all our solutions with respect to the
(+)-connection on the compact nilmanifold $M_3$ satisfy the
heterotic equations of motion.
\end{abstract}

\title[ Heterotic String Compactifications with non-zero fluxes and
constant dilaton] {Non-Kaehler Heterotic String Compactifications\\
with non-zero fluxes and constant dilaton}
\date{\today}

\author{Marisa Fern\'andez}
\address[Fern\'andez]{Universidad del Pa\'{\i}s Vasco\\
Facultad de Ciencia y Tecnolog\'{\i}a, Departamento de Matem\'aticas\\
Apartado 644, 48080 Bilbao\\ Spain} \email{marisa.fernandez@ehu.es}

\author{Stefan Ivanov}
\address[Ivanov]{University of Sofia "St. Kl. Ohridski"\\
Faculty of Mathematics and Informatics\\
Blvd. James Bourchier 5\\
1164 Sofia, Bulgaria} \email{ivanovsp@fmi.uni-sofia.bg}

\author{Luis Ugarte}
\address[Ugarte]{Departamento de Matem\'aticas\,-\,I.U.M.A.\\
Universidad de Zaragoza\\
Campus Plaza San Francisco\\
50009 Zaragoza, Spain} \email{ugarte@unizar.es}

\author{Raquel Villacampa}
\address[Villacampa]{Departamento de Matem\'aticas\,-\,I.U.M.A.\\
Universidad de Zaragoza\\
Campus Plaza San Francisco\\
50009 Zaragoza, Spain} \email{raquelvg@unizar.es}

\maketitle

\setcounter{tocdepth}{2} \tableofcontents

\section{Introduction. Field and Killing-spinor equations}

The bosonic fields of the ten-dimensional supergravity which arises
as low energy effective theory of the heterotic string are the
spacetime metric $g$, the NS three-form field strength $H$, the
dilaton $\phi$ and the gauge connection $A$ with curvature $F^A$.
The bosonic geometry is of the form $R^{1,9-d}\times M^d$ where the
bosonic fields are non-trivial only on $M^d$, $d\leq 8$. We consider
the two connections
\begin{equation*}
\nabla^{\pm}=\nabla^g \pm \frac12 H,
\end{equation*}
where $\nabla^g$ is the Levi-Civita connection of the Riemannian
metric $g$. Both connections preserve the metric, $\nabla^{\pm}g=0$ and
have totally skew-symmetric torsion $\pm H$, respectively.


The Green-Schwarz anomaly cancellation mechanism requires that the
three-form Bianchi identity receive an $\alpha'$ correction of the
form
\begin{equation}\label{acgen}
dH=\frac{\alpha'}48\pi^2(p_1(M^p)-p_1(E))=\frac{\alpha'}4
\Big(Tr(R\wedge R)-Tr(F^A\wedge F^A)\Big),
\end{equation}
where $p_1(M^p), p_1(E)$ are the first Pontrjagin forms of $M^p$
with respect to a connection $\nabla$ with curvature $R$ and the
vector bundle $E$ with connection $A$, respectively.

A class of heterotic-string backgrounds for which the Bianchi
identity of the three-form $H$ receives a correction of type
\eqref{acgen} are those with (2,0) world-volume supersymmetry. Such
models were considered in \cite{HuW}. The target-space geometry of
(2,0)-supersymmetric sigma models has been extensively investigated
in \cite{HuW,Str,HP1}. Recently, there is revived interest in these
models \cite{GKMW,CCDLMZ,GMPW,GMW,GPap} as string backgrounds and in
connection to heterotic-string compactifications with fluxes
\cite{Car1,BBDG,BBE,BBDP,y1,y2,y3,y4}.

In writing \eqref{acgen} there is a subtlety to the choice of
connection $\nabla$ on $M^p$ since anomalies can be cancelled
independently of the choice \cite{Hull}. Different connections
correspond to different regularization schemes in the
two-dimensional worldsheet non-linear sigma model. Hence the
background fields given for the particular choice of $\nabla$ must
be related to those for a different choice by a field redefinition
\cite{Sen}. Connections on $M^d$ proposed to investigate the anomaly
cancellation  \eqref{acgen} are $\nabla^g$ \cite{Str,GMW},
$\nabla^+$ \cite{CCDLMZ} and very recently \cite{DFG}, $\nabla^-$ \cite{Hull,Berg,Car1,GPap,II,KY}, Chern
connection $\nabla^c$ when $d=6$ \cite{Str,y1,y2,y3,y4}.

A heterotic geometry will preserve supersymmetry if and only if, in
10 dimensions, there exists at least one Majorana-Weyl spinor
$\epsilon$ such that the supersymmetry variations of the fermionic
fields vanish, i.e. the next Killing-spinor equations hold
\cite{Str}
\begin{gather} \nonumber
\delta_{\lambda}=\nabla_m\epsilon = \left(\nabla_m^g
+\frac{1}{4}H_{mnp}\Gamma^{np} \right)\epsilon=\nabla^+\epsilon=0
\\\label{sup1} \delta_{\Psi}=\left(\Gamma^m\partial_m\phi
-\frac{1}{12}H_{mnp}\Gamma^{mnp} \right)\epsilon=(d\phi-\frac12H)\cdot\epsilon=0 \\
\nonumber
\delta_{\xi}=F^A_{mn}\Gamma^{mn}\epsilon=F^A\cdot\epsilon=0,
\end{gather}
where  $\lambda, \Psi, \xi$ are the gravitino, the dilatino and the
gaugino, fields, respectively and $\cdot$ means Clifford action of
forms on spinors.

The bosonic part of the ten-dimensional supergravity action in the
string frame is (\cite{Berg}, $R=R^-$)
\begin{gather}\label{action}
S=\frac{1}{2k^2}\int
d^{10}x\sqrt{-g}e^{-2\phi}\Big[Scal^g+4(\nabla^g\phi)^2-\frac{1}{2}|H|^2
-\frac{\alpha'}4\Big(Tr |F^A|^2)-Tr |R|^2\Big)\Big].
\end{gather}

The string frame field equations (the equations of motion induced
from the action \eqref{action}) of the heterotic string up to
two-loops \cite{HT} in sigma model perturbation theory are (we use
the notations in \cite{GPap})
\begin{gather}\nonumber
Ric^g_{ij}-\frac14H_{imn}H_j^{mn}+2\nabla^g_i\nabla^g_j\phi-\frac{\alpha'}4
\Big[(F^A)_{imab}(F^A)_j^{mab}-R_{imnq}R_j^{mnq}\Big]=0;\\\label{mot}
\nabla^g_i(e^{-2\phi}H^i_{jk})=0;\\\nonumber
\nabla^+_i(e^{-2\phi}(F^A)^i_j)=0.
\end{gather}
The field equation of the dilaton $\phi$ is implied from the first
two equations above.

The  first compact torsional solutions for the heterotic/type I
string were obtained via duality from M-theory compactifications on
${\rm K3}\times {\rm K3}$ proposed in \cite{DRS}.  The metric was
first written down on the orientifold limit in \cite{DRS} and such
backgrounds have since been studied (see \cite{BBDG,BBE} and
references therein). The metric and the $H$-flux are derived by
applying a chain of supergravity dualities and the resulting
geometry in the heterotic theory is a $ \mathbb T^2$ bundle over a
K3.

Compact example solving \eqref{sup1}  and \eqref{acgen} with nonzero
field strength, constant dilaton and taking $R=R^+$, is constructed
in \cite{CCDLMZ} on the Iwasawa nilmanifold which is a $\mathbb T^2$
bundle over $\mathbb T^4$. However, it has been pointed out in
\cite{GMW} that this example is not a valid solution due to a sign
error in the torsional equation derived from the first two equations
in \eqref{sup1} which leads to the opposite sign in the left hand
side of \eqref{acgen}. Compact example of a balanced 6-manifold with constant dilaton 
non-trivial warped factor and  torsion generated by the Chern-Simons term only is presented
very recently in \cite{DFG}.

Compact examples in dimension six solving \eqref{sup1} and
\eqref{acgen} with non-zero flux $H$ and non-constant dilaton were
constructed by Li and Yau \cite{y1} for U(4) and U(5) principal
bundles taking $R=R^c$-the curvature of the Chern connection in
\eqref{acgen}. Non-Kaehler compact solutions of \eqref{sup1} and \eqref{acgen} on some torus bundles
over Calabi-Yau 4-manifold (K3 surfaces or complex torus) provided
in \cite{GP} are presented by Fu and Yau \cite{y2,y3} using the
Chern connection in \eqref{acgen}. It is confirmed in \cite{y4} that
the examples of torus bundles over the complex torus can not be
solutions to \eqref{sup1} and \eqref{acgen} taking with respect to
the curvature of the Chern connection $ R=R^c$ with $\alpha'>0$
while some torus bundles over K3 surfaces are valid solutions.

It is known \cite{Bwit,GMPW} (\cite{GPap} for dimension 6) that the
equations of motion of type I supergravity are automatically
satisfied with $R=0$ if one imposes, in addition to the preserving
supersymmetry equations \eqref{sup1}, the three-form Bianchi
identity \eqref{acgen} taking with respect to a flat connection on
$TM$, $R=0$.

According to no-go (vanishing) theorems  (a consequence of the
equations of motion \cite{FGW,Bwit}; a consequence of the
supersymmetry \cite{IP1,IP2} for SU($n$)-case and \cite{GMW} for the
general case) there are no compact solutions with non-zero flux and
non-constant dilaton satisfying simultaneously the supersymmetry
equations \eqref{sup1} and the three-form  Bianchi identity
\eqref{acgen} with  $Tr (R\wedge R)=0$.

However, in the presence of a curvature term $R$ the solution of the
supersymmetry equations \eqref{sup1} and the anomaly cancellation
condition \eqref{acgen} obey the second and the third equations of
motion but do not always satisfy the gravitino equation of motion
(the first equation in \eqref{mot}). If $R$ is an SU($3$)-instanton
then \eqref{sup1} and \eqref{acgen} imply \eqref{mot}. This can be
seen from the considerations in the Appendix of \cite{GMPW}. 
We give a quadratic expression for $R$ which is necessary and sufficient condition in order that 
\eqref{sup1} and \eqref{acgen} imply \eqref{mot}  based on the properties of the
special geometric structure induced from the first two equations in
\eqref{sup1}. More precisely, we prove in Section \ref{prof} the
following
\begin{thrm}\label{thac}
Let $(M,J,g,F^A,R)$ be a conformally balanced Hermitian manifold
with Kaehler form $F$ which solves the heterotic Killing spinor
equations \eqref{sup1} and the anomaly cancellation \eqref{acgen}.
\begin{description}
\item[a)]
The gravitino equation of motion  is a consequence of the heterotic
Killing spinor equations \eqref{sup1} and  the anomaly cancellation
\eqref{acgen}  if and only if the next identity holds
\begin{equation}\label{supmot}
\frac1{2}\Big[R_{msab}R_{pqab}+R_{mpab}R_{qsab}+R_{mqab}R_{spab}\Big]F^{pq}J^s_n
=R_{mpqr}R_n^{pqr}.
\end{equation}
\begin{itemize}
\item If $R$ is $(1,1)$-form, $J^s_mJ^p_nR_{spab}= R_{mnab}$ then \eqref{supmot} is equivalent to
\begin{equation}\label{supmot6}
R_{mjab}R_{klab}F^{kl} =0.
\end{equation}
In particular, the gravitino equation of motion with respect to
either the Chern connection or the (--)-connection is a consequence
of the heterotic Killing spinor equations \eqref{sup1} and  the
anomaly cancellation \eqref{acgen} if and only if \eqref{supmot6}
holds.
\item
If $R$ is an {\rm SU($3$)}-instanton then \eqref{supmot} holds.
\end{itemize}
\item[b)]
If $R^-$ is an {\rm SU($3$)}-instanton, $Hol(\nabla^+)\subset {\frak
s}{\frak u}(3)$ and the manifold is compact then the flux $H$
vanishes, the dilaton is constant and the manifold is a Calabi-Yau
space.
\end{description}
\end{thrm}
As a consequence of Theorem \ref{thac}, considering solutions involving Chern connection, one may 
study  stability of the tangent bundle.

The main goal of this paper is to construct explicit compact valid
solutions with non-zero field strength, non-flat instanton and
constant dilaton to the heterotic equations of motion \eqref{mot} in
dimension six. We present compact nilmanifolds in dimension six
satisfying the heterotic supersymmetry equations \eqref{sup1} with
non-zero flux $H\not=0$, non-flat instanton $F^A\not=0$ and constant
dilaton obeying the three-form Bianchi identity \eqref{acgen} with
curvature term $R=R^g$, $R=R^+$ or $R=R^c$. Some of them are torus
bundles over the complex torus but this does not violate the
non-existence result in \cite{y4} since we use different curvature
term ($R^g$ or $R^+$) in \eqref{acgen}. In particular, we present a
valid solution on the Iwasawa manifold but with respect to a
non-standard complex structure. We find compact valid solutions to
\eqref{sup1} with non-zero flux, non-flat instanton  and constant
dilaton satisfying anomaly cancellation condition \eqref{acgen}
using the curvature $R^c$ of the Chern connection on an $S^1$ bundle
over a 5-manifold which is a $\mathbb T^2$ bundle over $\mathbb
T^3$. All manifolds do not admit any Kaehler metric and seem to be
the first explicit compact valid supersymmetric heterotic solutions
to \eqref{sup1}  and \eqref{acgen}  with non-zero flux, non-flat
instanton and constant dilaton in dimension six.

However, because of Theorem \ref{thac}, the gravitino equation of
motion is not satisfied in most cases. Only the  solutions
constructed in Theorem~\ref{solutions-h3-1} b),
Theorem~\ref{solutions-h3-2} b) on the compact nilmanifold
$M_3=\Gamma\backslash H(2,1)\times S^1$, where $H(2,1)$ is the
5-dimensional Heisenberg group and $\Gamma$ is a lattice, solve in
addition the heterotic equations of motion \eqref{mot} with non-zero
fluxes and constant dilaton. It seems that these are the first
compact supersymmetric solutions to the heterotic equations of
motion with non-zero flux $H\not=0$, non-flat instanton $F^A\not=0$
and constant dilaton in dimension six.

Our convention for the curvature is given in Section \ref{pre}.

\section{The supersymmetry equations in dimensions 6}

Necessary and sufficient conditions to have a solution to the system
of dilatino and gravitino equations in dimension 6 were derived by
Strominger in \cite{Str} involving the notion of SU($n$)-structure
and then studied by many authors
\cite{GKMW,GMPW,GMW,CCDLMZ,Car1,II,BBDG,BBE,GPap,y1,y2,y3,y4}.

\subsection{SU(3)-structures in $d=6$}
Let $(M,J,g)$ be an almost Hermitian 6-manifold with Riemannian
metric $g$ and almost complex structure $J$, i.e. $(J,g)$ define an
$U(3)$-structure. The Nijenhuis tensor $N$, the Kaehler form $F$ and
the Lee form $\theta$ are defined by
\begin{equation*}\label{cy1}
N(\cdot,\cdot)=[J\cdot,J\cdot]-[\cdot,\cdot]-J[J\cdot,\cdot]-J[\cdot,J\cdot],
\quad F(\cdot,\cdot)=g(\cdot,J\cdot), \quad \theta(\cdot)=\delta
F(J\cdot),
\end{equation*}
respectively, where $*$ is the Hodge operator and $\delta$ is the
co-differential, $\delta=-*d*$.

An SU($3$)-structure is determined by an additional non-degenerate
(3,0)-form $\Psi=\Psi^+ + i\, \Psi^-$, or equivalently by a
non-trivial spinor.
The subgroup of $SO(6)$ fixing the forms $F$ and  $\Psi$
simultaneously is SU($3$).  The Lie algebra of SU($3$) is denoted
${\frak s}{\frak u}(3)$.

The failure of the holonomy group of the Levi-Civita connection to
reduce to SU($3$) can be measured by the intrinsic torsion $\tau$,
which is identified with $\nabla^g F$ or $\nabla^g J$ and can be
decomposed into five classes \cite{CS}, $\tau \in W_1\oplus \dots
\oplus W_5$. The intrinsic torsion of an $U(n)$-structure belongs to
the first four components described by Gray-Hervella \cite{GrH}. The
five components of an SU($3$)-structure are first described by
Chiossi-Salamon \cite{CS} (for interpretation in physics see
\cite{CCDLMZ}) and are determined by $dF,d\Psi^+,d\Psi^-$ as well as
by $dF$ and $N$. The Hermitian manifolds belong to $W_3\oplus W_4$.
In the paper we are interested in the class $W_3$ of balanced
Hermitian manifolds \cite{Mi} characterized by the conditions $N=0,
\theta=0$ or, equivalently, $N=0$, $d*F=0$.

Necessary conditions to solve the gravitino equation are given in
\cite{FI}. The presence of a parallel spinor in dimension 6 leads
firstly to the reduction to $U(3)$, i.e. the existence of an almost
Hermitian structure, secondly to the existence of a linear
connection preserving the almost Hermitian structure with torsion
3-form and thirdly to the reduction of the holonomy group of the
torsion connection to SU($3$), i.e. its Ricci 2-form has to be
identically zero. It is shown in \cite{FI} that there exists a
unique linear connection preserving an almost Hermitian structure
having totally skew-symmetric torsion  if and only if the Nijenhuis
tensor is a 3-form, i.e. the intrinsic torsion $\tau \in W_1\oplus
W_3\oplus W_4$. The torsion connection $\nabla^+$ with torsion $T$
is determined by
\begin{equation*}\label{cy2}
\nabla^+ = \nabla^g +\frac{1}{2}T, \qquad
T=JdF+N=-dF(J\cdot,J\cdot,J\cdot)+N.
\end{equation*}

Necessary and sufficient conditions to solve the gravitino equation
in dimension 6 are given in \cite{II}. Namely, there exists a unique
linear connection with torsion 3-form which preserves the almost
Hermitian structure whose holonomy is contained in SU($3$) if and
only if the first Chern class vanishes, $c_1(M,J)=0$ and the
SU($3$)-structure $(M,g,F,\Psi^+,\Psi^-)$ satisfies the differential
equations \cite{II}
\begin{equation}\label{cycon}
d\Psi^+=\theta\wedge\Psi^+ -\frac{1}{4}(N,\Psi^+)*F, \qquad
d\Psi^-=\theta\wedge\Psi^- -\frac{1}{4}(N,\Psi^-)*F.
\end{equation}

The torsion $T$ is given by $T=-*dF +*(\theta\wedge
F)+\frac{1}{4}(N,\Psi^+)\Psi^+ + \frac{1}{4}(N,\Psi^-)\Psi^-.$

Necessary and sufficient conditions to solve the gravitino and
dilatino equations are presented in \cite{Str}. The dilatino
equation forces the almost complex structure to be integrable
($N=0$) and the Lee form to be closed (for applications in physics
the Lee form has to be exact) determined by the dilaton due to
$\theta=2d\phi$ \cite{Str}. The three-form field strength
$H=T=-dF(J\cdot,J\cdot,J\cdot)=-*dF+*(2d\phi\wedge F)$. Solutions
with constant dilaton are those with zero Lee form, $dF^{n-1}=0$,
i.e. balanced Hermitian manifolds.

When the almost complex structure is integrable, $N=0$, the torsion
connection $\nabla^+$ is also known as the Bismut connection (we
shall call it Bismut-Strominger (B-S) connection) and was used by
Bismut to prove local index theorem for the Dolbeault operator on
non-Kaehler Hermitian manifolds \cite{Bis}. This formula was
recently applied in string theory \cite{BBE}. Vanishing theorems for
the Dolbeault cohomology on compact non-Kaehler Hermitian manifolds
were found in terms of the B-S connection \cite{AI,IP2,IP1}.

In addition to these equations, the vanishing of the gaugino
variation requires the non-zero 2-form $F^A$ to be of instanton type
(\cite{CDev,Str,GMW}). A Donaldson-Uhlenbeck-Yau SU($3$)-instanton
i.e. the gauge field $A$ is a connection on a holomorphic vector
bundle with curvature 2-form $F^A\in {\frak s}{\frak u}(3)$.
The SU($3$)-instanton condition can be written in local holomorphic
coordinates in the form \cite{CDev,Str}
\begin{equation*}\label{6inst}
F^A_{\alpha\beta}=F^A_{\bar{\alpha}\bar{\beta}}=0,\quad\
F^A_{\alpha\bar\beta}F^{\alpha\bar\beta}=0.
\end{equation*}

\subsection{Proof of Theorem \ref{thac} }\label{prof}

A consequence of the gravitino and dilatino Killing equation is the
expression of the Ricci tensor $Ric^+_{mn}=R^+_{imnj}g^{ij}$ of the
(+)-connection established in \cite{IP1}, Proposition 3.1:
\begin{equation}\label{ric+ff}
Ric^+_{mn}=-2\nabla^+_md\phi_n -\frac14dT_{mspq}J_n^sF^{pq}
=-2\nabla^g_md\phi_n+d\phi_sT^s_{mn} -\frac14dT_{mspq}J_n^sF^{pq}.
\end{equation}
The four-form $dT=dJdF$ is a (2,2)-form with respect to the complex
structure $J$. Therefore, the last term in \eqref{ric+ff} is
symmetric.

On the other hand,  the Ricci tensors of $\nabla^g$ and $\nabla^+$
are connected by (see e.g. \cite{FI})
\begin{gather}\label{ricg+}
Ric^g_{mn}=Ric^+_{mn}+\frac14T_{mpq}T_n^{pq}-\frac12\nabla^+_sT^s_{mn},
\qquad
Ric^+_{mn}-Ric^+_{nm}=\nabla^+_sT^s_{mn}=\nabla^g_sT^s_{mn},\\\label{mo}
Ric^g_{mn}=\frac12(Ric^+_{mn}+Ric^+_{nm})+\frac14T_{mpq}T_n^{pq}.
\end{gather}
Substitute \eqref{ric+ff} into \eqref{mo}, insert the result into
the first equation of \eqref{mot} and use the anomaly cancellation
\eqref{acgen} to conclude \eqref{supmot}.  If $R$ is a (1,1)-form
then \eqref{supmot6} is a consequence of \eqref{supmot}. It is well
known that the curvature of the Chern connection $R^c$ is always a
(1,1)-form. When $Hol(\nabla^+)\subset {\frak s}{\frak u}(3)$ the
curvature  $R^-$ of the (--)-connection is also an (1,1)-form. This
follows from the well known identity
\begin{equation}\label{r1}
dT_{ijkl}=2R^+_{ijkl}-2R^-_{klij}
\end{equation}
 and the fact that $dT$ is a (2,2)-form. This completes the proof of  a).

The proof of b) is essentially contained in \cite{IP1,IP2}. Indeed,
if $Hol(\nabla^+)\subset {\frak s}{\frak u}(3)$  and $R^-$ is an
SU($3$)-instanton, \eqref{r1} yields $dT_{ispq}F^{pq}=0$, i.e. the
manifold is «almost strong« in the terminology of  \cite{IP1}. Then
Corollary 4.2 a) in \cite{IP1} asserts that there are no holomorphic
(3,0) forms which contradicts the result in  \cite{Str} except
$T=d\phi=0$. This completes the proof of Theorem \ref{thac}.

\subsection{Heterotic supersymmetry with constant
dilaton}\label{cdil}

We look for a compact Hermitian 6-manifold $(M,J,g)$ which satisfies
the following conditions
\begin{enumerate}
\item[a).] Gravitino equation: $Hol(\nabla^+)\subset {\frak
s}{\frak u}(3)$, i.e.
\begin{equation}\label{1}\sum_{i=1}^6(\Omega^+)^{E_i}_{JE_i}=0,
\end{equation}
where $\{E_1,\ldots,E_6\}$ is an orthonormal basis on $M$.
\item[b).] Dilatino equation with constant dilaton: the Lee form
$\theta=2d\phi=0$, i.e. $(M,J,g)$ is a balanced manifold.
\item[c).] Gaugino equation: look for a Hermitian vector bundle
$E$ of rank $r$ over $M$ equiped with an SU($3$)-instanton, i.e. a
connection $A$ with curvature 2-form $\Omega^A$ satisfying
\begin{equation}\label{2}
(\Omega^A)^i_j(JE_k,JE_l)=(\Omega^A)^i_j(E_k,E_l),\qquad
\sum_{k=1}^6(\Omega^A)^i_j(E_k,JE_k)=0.
\end{equation}
\item[d).] Anomaly cancellation condition:
\begin{equation}\label{ac}
dH=dT=\frac{\alpha'}48\pi^2\Big(p_1(M)-p_1(A)\Big), \qquad
\alpha'>0.
\end{equation}
\end{enumerate}

\section{General preliminaries}\label{pre}

For a linear connection $\nabla$, the connection 1-forms
$\omega^i_j$ with respect to a fixed basis $E_1,\dots,E_6$ are
$$
\omega^i_j(E_k) = g(\nabla_{E_k}E_j,E_i)
$$
since we write $\nabla_X E_j = \omega^1_j(X)\, E_1 +\cdots+
\omega^6_j(X)\, E_6$.

The curvature 2-forms $\Omega^i_j$ of  $\nabla$ are given in terms
of the connection 1-forms $\omega^i_j$ by
$$
\Omega^i_j  = d \omega^i_j + \omega^i_k\wedge\omega^k_j, \quad
 \Omega_{ji} = d \omega_{ji} +
\omega_{ki}\wedge\omega_{jk}, \quad R^l_{ijk}=\Omega^l_k(E_i,E_j),
\quad R_{ijkl}=R^s_{ijk}g_{ls}.
$$
and the first Pontrjagin class is represented by the 4-form
$$
p_1(\nabla)={1\over 8\pi^2} \sum_{1\leq i<j\leq 6}
\Omega^i_j\wedge\Omega^i_j.
$$

Let $(M,J,g)$ be a 6-dimensional Hermitian manifold. Consider the
connections with torsion $\nabla^{\pm}$ given by $\nabla^{\pm}=\nabla^g \pm
\frac12 T$ with torsion $T$ given by
\begin{equation}\label{tor}
 T=JdF= -*dF.
\end{equation}
Notice that $\nabla^+$ is precisely the B-S connection of the
Hermitian structure.

The Chern connection $\nabla^c$ is defined by,
$$
\nabla^c= \nabla^g + \frac12 C, \qquad C(.,.,.)=dF(J.,.,.).
$$
Observe that the tensor field $C$ satisfies that
$C(X,\cdot,\cdot)=(JX\lrcorner dF)(\cdot,\cdot)$ is a 2-form on $M$.

Let us suppose that $(J,g)$ is a left invariant Hermitian structure
on a 6-dimensional Lie group $M$ and let $\{e^1,\ldots,e^6\}$ be an
orthonormal basis of left invariant 1-forms, that is, $g=e^1\otimes
e^1 + \cdots +e^6\otimes e^6$. Let
$$
d\, e^k = \sum_{1\leq i<j \leq 6} a_{ij}^k \, e^i\wedge
e^j,\quad\quad k=1,\ldots,6,
$$
be the structure equations in the basis $\{e^k\}$.

Let us denote by $\{E_1,\ldots,E_6\}$ the dual basis. Since $d
e^k(E_i,E_j)= -e^k([E_i,E_j])$, we have that the Levi-Civita
connection 1-forms $(\omega^g)^i_j$ are
\begin{equation}\label{lc}
(\omega^g)^i_j(E_k) = -\frac12(g(E_i,[E_j,E_k]) - g(E_k,[E_i,E_j]) +
g(E_j,[E_k,E_i]))=\frac12(a^i_{jk}-a^k_{ij}+a^j_{ki}).
\end{equation}

The connection 1-forms $(\omega^\pm)^i_j$ for the connections with
torsion $\nabla^\pm$ are given by
\begin{equation}\label{pm}
(\omega^\pm)^i_j(E_k)=(\omega^g)^i_j(E_k) + \frac12
(T^{\mp})^i_j(E_k), \qquad (T^{\pm})^i_j(E_k)=T^{\pm}(E_i,E_j,E_k)=\mp
dF(JE_i,JE_j,JE_k).
\end{equation}

The connection 1-forms $(\omega^c)^i_j$ for the Chern connection
$\nabla^c$ are determined by
\begin{equation}\label{c}
(\omega^c)^i_j(E_k)=(\omega^g)^i_j(E_k) + \frac12 C^i_j(E_k), \qquad
C^i_j(E_k) = dF(JE_k,E_i,E_j).
\end{equation}

We shall focus on six-dimensional nilmanifolds $M=\Gamma\backslash
G$ endowed with an invariant (integrable almost) complex structure
$J$. According to Proposition~6.1 in \cite{FPS}, for invariant
Hermitian metrics on compact nilmanifolds the balanced condition is
equivalent to  $Hol(\nabla^+)\subset {\frak s}{\frak u}(3)$. The
equivalence of the conditions a) and b) in subsection~\ref{cdil} can
also be derived from \eqref{cycon} and the fact, established in
\cite{FPS}, that for any invariant Hermitian structure on a
nilmanifold the (3,0)-form $\Psi=\Psi^+ + i\, \Psi^-$ is closed.

\subsection{Six-dimensional balanced Hermitian nilmanifolds}
Next we review the main results given in~\cite{U} concerning
balanced $J$-Hermitian metrics on $M$ in order to apply them to the
construction of solutions to the equations \eqref{1}-\eqref{ac}
above. First of all, if $(M,J)$ admits a balanced $J$-Hermitian
metric (not necessarily invariant) then the Lie algebra $\frg$ of
$G$ is isomorphic to $\frh_{1},\ldots,\frh_6$ or $\frh_{19}^-$,
where $\frh_{1} = (0,0,0,0,0,0)$ is the abelian Lie algebra and
$$
\begin{array}{rcl}
\frh_{2} &\!\!=\!\!& (0,0,0,0,12,34),\\[2pt]
\frh_{3} &\!\!=\!\!& (0,0,0,0,0,12+34),\\[2pt]
\frh_{4} &\!\!=\!\!& (0,0,0,0,12,14+23),
\end{array}
\quad\quad\quad
\begin{array}{rcl}
\frh_{5} &\!\!=\!\!& (0,0,0,0,13+42,14+23),\\[2pt]
\frh_{6} &\!\!=\!\!& (0,0,0,0,12,13),\\[2pt]
\frh^-_{19} &\!\!=\!\!& (0,0,0,12,23,14-35).
\end{array}
$$

Here $\frh_5$ is the Lie algebra underlying the Iwasawa manifold.
For the canonical complex structure $J_0$ on $\frh_5$ there exists a
complex basis $\{\omega^j\}_{j=1}^3$ of 1-forms of type (1,0)
satisfying $d\omega^1= d\omega^2=0$ and $d\omega^3= \omega^{12}$.

Since the Lie algebras $\frh_2,\ldots,\frh_6$ are 2-step nilpotent,
for any complex structure $J$ ($\not= J_0$ for $\frh_5$) there is a
basis $\{\omega^j\}_{j=1}^3$ of (1,0)-forms such that
\begin{equation}\label{nilp-struct}
d \omega^1=d\omega^2=0,\quad\quad d\omega^3=\rho\, \omega^{12} +
\omega^{1\bar{1}} + B\,\omega^{1\bar{2}} + D\,\omega^{2\bar{2}},
\end{equation}
where $B,D\in \mathbb{C}$, and $\rho= 0,1$. In particular, $J$ is a
nilpotent complex structure on $\frh_2,\ldots,\frh_6$ in the sense
\cite{CFGU}. Recall that a complex structure $J$ on a
$2n$-dimensional nilpotent Lie algebra $\frg$ is called {\em
nilpotent} if there is a basis $\{\omega^j\}_{j=1}^n$ of (1,0)-forms
satisfying $d\omega^1 = 0$ and
$$
d\omega^j \in {\bigwedge} ^{2} (\omega^1, \ldots, \omega^{j-1},
\omega^{\overline{1}}, \ldots, \omega^{\overline{j-1}}),
$$
for $j=2, \cdots,n$.

Any complex structure on the Lie algebra $\frh^-_{19}$ is not
nilpotent and there is a (1,0)-basis $\{\omega^j\}_{j=1}^3$
satisfying
\begin{equation}\label{nonilp-struct}
d\omega^1 = 0,\quad d\omega^2 = E\, \omega^{13} +
\omega^{1\bar{3}},\quad d\omega^3 = C\, \omega^{1\bar{1}} + ia\,
\omega^{1\bar{2}} - ia\bar{E}\, \omega^{2\bar{1}},
\end{equation}
where $E\in \mathbb{C}$ with $|E|=1$, $\bar{C}=CE$ and
$a\in\mathbb{R}-\{0\}$.

Now, the fundamental form $F$ of any invariant $J$-Hermitian
structure is given in terms of the basis $\{\omega^j\}_{j=1}^3$ by
\begin{equation}\label{2forma}
2\, F=i (r^2\omega^{1\bar{1}} + s^2\omega^{2\bar{2}} +
t^2\omega^{3\bar{3}})
+u\,\omega^{1\bar{2}}-\bar{u}\,\omega^{2\bar{1}} +
v\,\omega^{2\bar{3}}-\bar{v}\,\omega^{3\bar{2}}+
z\,\omega^{1\bar{3}}-\bar{z}\,\omega^{3\bar{1}},
\end{equation}
where $r,s,t\in\mathbb{R}-\{0\}$ and $u,v,z\in\mathbb{C}$ must
satisfy those restrictions coming from the positive definiteness of
the associated metric $g(X,Y)=-F(X,JY)$. The following result gives
necessary and suficient conditions, in terms of the different
coefficients involved, in order the Hermitian structure to be
balanced.

\begin{prop} \label{balanced} \cite{U}
In the notation above, we have:
\begin{enumerate}
\item[{\rm ({\it i})}] If $J$ is a nonnilpotent complex structure
defined by~$(\ref{nonilp-struct})$, then $(J,F)$ is balanced if and
only if
$$
z=-iuv/s^2 \quad\quad \mbox{and} \quad\quad
Cs^2+a\bar{E}u+a\bar{u}=0.
$$
\item[{\rm ({\it ii})}] If $J$ is a nilpotent complex
structure defined by~$(\ref{nilp-struct})$, then $(J,F)$ is balanced
if and only if
$$
s^2t^2-|v|^2 + D(r^2t^2-|z|^2) = B(it^2\bar{u} - v\bar{z}).
$$
\end{enumerate}
\end{prop}

\section{The Iwasawa manifold revisited}\label{Iwasawa}

Apart from the abelian Lie algebra, $\frh_5$ is the only
6-dimensional nilpotent Lie algebra which can be given a {\it
complex} Lie algebra structure. The corresponding complex
parallelizable nilmanifold is the well-known Iwasawa manifold. This
manifold is studied in~\cite{CCDLMZ}; however, as it is pointed out
in the introduction, this example is not a valid solution due to a
sign error in the torsional equation. Here we show that there are no
valid solutions on the Iwasawa manifold, a fact which leads us to
study more general complex nilmanifolds in the subsequent sections.

The standard complex structure $J_0$ on $\frh_5$ is defined by the
following complex structure equations:
$$
d\omega^1=d\omega^2=0,\quad\quad d\omega^3 = \omega^{12}.
$$
For any $t\not=0$, let us consider $F$ given by
$$
F=\frac{i}{2} (\omega^{1\bar{1}} + \omega^{2\bar{2}} +
t^2\,\omega^{3\bar{3}}).
$$
It is easy to see that the Hermitian structure $(J_0,F)$ is balanced
for any value of the parameter.

Notice that the Iwasawa manifold is a $\mathbb{T}^2$ bundle over
$\mathbb{T}^4$, where the parameter $t$ scales the fiber.

From a real point of view, let us consider the real basis of 1-forms
$\{e^1,\ldots,e^6\}$ given by
$$
e^1+i\, e^2=\omega^1,\quad e^3+i\, e^4=\omega^2,\quad e^5+i\,
e^6=t\, \omega^3.
$$
Now, in terms of this basis, we have that the structure equations
are
\begin{equation}\label{family-Iwasawa}
\left\{
  \begin{aligned}
  &d e^1= d e^2= d e^3= d e^4= 0, \\
  &d e^5= t\, e^{13}- t\, e^{24},\\
  &d e^6= t\, e^{14}+ t\, e^{23},
  \end{aligned}
\right.
\end{equation}
the complex structure $J_0$ is given by $J_0 e^1=-e^2,  J_0
e^3=-e^4, J_0 e^5=-e^6,$ the $J_0$-Hermitian metric $g=e^1\otimes
e^1+\cdots + e^6\otimes e^6$ has the associated fundamental form
$F=e^{12}+e^{34}+e^{56}.$ The structure equations
\eqref{family-Iwasawa} give $dF= t\, e^{136}- t\, e^{145}- t\,
e^{235}- t\, e^{246}$. Apply \eqref{tor} to verify that the torsion
$T$ of $\nabla^+$ satisfies
$$
T= -t\, e^{135}- t\, e^{146}- t\, e^{236}+ t\, e^{245}, \qquad dT-4 t^2 e^{1234}.$$

All the curvature forms $(\Omega^c)^i_j$ of the Chern connection
vanish. In view of \eqref{lc} and \eqref{pm}, the non-zero curvature
forms $(\Omega^g)^i_j$ and $(\Omega^\pm)^i_j$ for the Levi-Civita
connection and the connections $\nabla^\pm$ are given by:

$$
\begin{array}{l}
(\Omega^g)^1_2={t^2\over 2}(e^{34}-e^{56}),\quad
(\Omega^g)^1_3=-{t^2\over 4}(3e^{13}-e^{24}),\quad
(\Omega^g)^1_4=-{t^2\over 4}(3e^{14}+e^{23}),\\[5pt]
(\Omega^g)^1_5=-(\Omega^g)^2_6={t^2\over 4}(e^{15}-e^{26}), \quad
(\Omega^g)^1_6=(\Omega^g)^2_5={t^2\over 4}(e^{16}+e^{25}), \quad
(\Omega^g)^2_3=-{t^2\over 4}(e^{14}+3e^{23}),\\[5pt]
(\Omega^g)^2_4={t^2\over 4}(e^{13}-3e^{24}),\quad
(\Omega^g)^3_4={t^2\over 2}(e^{12}-e^{56}), \quad
(\Omega^g)^3_5=-(\Omega^g)^4_6={t^2\over 4}(e^{35}-e^{46}), \\[5pt]
(\Omega^g)^3_6=(\Omega^g)^4_5={t^2\over 4}(e^{36}+e^{45}),\quad
(\Omega^g)^5_6=-{t^2\over 2}(e^{12}+e^{34});
\end{array}
$$

$$
\begin{array}{l}
(\Omega^+)^1_2=2t^2 e^{34},\quad
(\Omega^+)^1_3=(\Omega^+)^2_4=-t^2(e^{13}+e^{24}),\quad
(\Omega^+)^1_4=-(\Omega^+)^2_3=-t^2(e^{14}-e^{23}),\\[5pt]
(\Omega^+)^3_4=2t^2e^{12}, \quad
(\Omega^+)^5_6=-2t^2(e^{12}+e^{34});
\end{array}
$$

$$
\begin{array}{l}
(\Omega^-)^1_2=(\Omega^-)^3_4=-2t^2 e^{56},\quad
(\Omega^-)^1_3=-(\Omega^-)^2_4=-t^2(e^{13}-e^{24}),\quad
(\Omega^-)^1_4=(\Omega^-)^2_3=-t^2(e^{14}+e^{23}).
\end{array}
$$

Clearly $Hol(\nabla^+)\subset {\frak s}{\frak u}(3)$ and the
Pontrjagin classes of the four connections are then represented by
\begin{equation}\label{p-J0-Iwasawa}
p_1(\nabla^g)={t^4\over 4\pi^2} e^{1234}, \quad
p_1(\nabla^+)=0,\quad p_1(\nabla^-)={t^4\over \pi^2} e^{1234},\quad
p_1(\nabla^c)=0.
\end{equation}

\subsection{Cardoso et al. abelian instanton}\label{card}
Cardoso et al. consider in~\cite{CCDLMZ} an abelian field strength
configuration with (1,1)-form
$$
\mathcal{F} = if\, dz_1\wedge d\bar{z}_1 - if\, dz_2\wedge
d\bar{z}_2 + {\rm e}^{i\gamma} \sqrt{\frac14 - f^2}\, dz_1\wedge
d\bar{z}_2 - {\rm e}^{-i\gamma} \sqrt{\frac14 - f^2}\, dz_2\wedge
d\bar{z}_1,
$$
where the function $f$ satisfies
$$
i\partial_{z_2} f + \partial_{z_1}\!\!\left({\rm e}^{-i\gamma}
\sqrt{\frac14 - f^2} \right)=0, \quad\quad\  i\partial_{z_1} f +
\partial_{z_2}\!\!\left({\rm e}^{i\gamma} \sqrt{\frac14 - f^2}
\right)=0.
$$
Under these conditions one gets
$$
Tr\ F\wedge F= \mathcal{F}\wedge \mathcal{F} = -\frac12 dz_1\wedge
dz_2\wedge d\bar{z}_1 \wedge d\bar{z}_2.
$$

Here $dz_1$ and $dz_2$ denote the (2,0)-forms at the level of the
Lie group, which descend to the forms $\omega^1$ and $\omega^2$ on
the compact nilmanifold. Therefore, on the  Iwasawa manifold we have
$$
Tr\ F\wedge F = \frac12\, \omega^{1\bar{1}} \wedge \omega^{2\bar{2}}
= -2\, e^{1234}.
$$

Now, taking $A$ as one of these abelian instantons we have that
\begin{equation}\label{cardoso-ex}
dT= -4 t^2 e^{1234} = -16\pi^2t^2 (p_1(\nabla^+)-p_1(A)),
\end{equation}
which is not a valid solution for any $t$ (see~\cite{GMW} for
details). Moreover, the whole space of complex structures compatible
with the canonical metric  obtained when $t=1$ in
\eqref{family-Iwasawa} is studied in~\cite{CCDLMZ} where the authors
proved that the behavior is the same as in \eqref{cardoso-ex}.

\begin{rmrk}\label{no-solution-J_0}
It is not difficult to prove that any $J_0$-Hermitian metric is
equivalent to one in the 1-parameter family given above. Since $dT-4t^2 e^{1234}$, in view of \eqref{p-J0-Iwasawa} there is no way to
find a satisfactory solution with $J_0$ as the underlying complex
structure. In fact, it is proved in \cite{y4} that torus bundles
over the complex torus can not be solutions to \eqref{sup1} and
\eqref{acgen} taking with respect to the curvature of the Chern
connection $R=R^c$ with $\alpha'>0$, which implies that $p_1(A)$
cannot be a positive multiple of $e^{1234}$ for any {\rm
SU($3$)}-instanton $A$ on the Iwasawa manifold, so \eqref{acgen}
cannot be satisfied neither for $R=R^g$ nor for $R=R^\pm$.
\end{rmrk}

Therefore, in order to find solutions we need to consider other
compact nilmanifolds or metrics and/or complex structures different
from the canonical ones on the nilmanifold underlying the Iwasawa
manifold. In the following sections we show many explicit solutions.

\section{A family of balanced Hermitian structures on the Lie algebra
$\frh_3$}\label{bal-family-h3}

In this section we construct explicit solutions on compact
nilmanifold corresponding to the Lie algebra $\frh_3$. First we
recall~\cite{U} that, up to equivalence, there exist two complex
structures $J^{\pm}$ on $\frh_3$, namely
$$
J^{\pm} \colon\ \ d\omega^1=d\omega^2=0, \quad d\omega^3 = \omega^{1\bar
1} \pm \omega^{2\bar 2},
$$
but only $J^-$ admits compatible balanced structures. Notice that
the balanced condition for $J^-$ given in
Proposition~\ref{balanced}~(ii) reduces to
$$
(r^2-s^2)t^2 = |z|^2 -|v|^2.
$$

For any $t\not=0$, let us consider the balanced structure $F$ given
by
$$
F=\frac{i}{2} (\omega^{1\bar{1}} + \omega^{2\bar{2}} +
t^2\,\omega^{3\bar{3}}),
$$
which corresponds to $r=s=1$ and $u=v=z=0$.

From a real point of view, let us consider the basis of 1-forms
$\{e^1,\ldots,e^6\}$ given by
$$
e^1+i\, e^2=\omega^1,\quad e^3+i\, e^4=\omega^2,\quad e^5+i\,
e^6=t\, \omega^3.
$$
Now, in terms of this basis, we have the structure equations
\begin{equation}\label{family-h3}
\left\{
  \begin{aligned}
  &d e^1= d e^2= d e^3= d e^4= d e^5 =0, \\
  &d e^6= -2t\, e^{12} + 2t\, e^{34},
  \end{aligned}
\right.
\end{equation}
and the complex structure $J=J^-$ is given by $ J e^1=-e^2, J
e^3=-e^4, J e^5=-e^6.$ The balanced $J$-Hermitian metric $
g=e^1\otimes e^1+\cdots + e^6\otimes e^6$ has the associated
fundamental form $F=e^{12}+e^{34}+e^{56}$. The structure equations
\eqref{family-h3} yield $dF=2t (e^{12} - e^{34})e^5$. For the
torsion $T$ of $\nabla^+$ we calculate using \eqref{tor}, \eqref{lc}
and \eqref{pm} that
$$
T= -2t(e^{12}-e^{34})e^6, \qquad dT=-8t^2 e^{1234}, \qquad \nabla^+
\, T=0.
$$

A direct calculation applying \eqref{lc} and \eqref{pm} shows that
the non-zero curvature forms $(\Omega^g)^i_j$ of the Levi-Civita
connection $\nabla^g$ are given by
$$
\begin{array}{l}
(\Omega^g)^1_2=-t^2 (3\, e^{12} - 2\, e^{34}),\quad
(\Omega^g)^1_3=t^2 e^{24},\quad (\Omega^g)^1_4=- t^2 e^{23},\quad
(\Omega^g)^1_6= t^2 e^{16},\quad
(\Omega^g)^2_3= -t^2 e^{14},  \\[10pt]
(\Omega^g)^2_4= t^2 e^{13}, \quad (\Omega^g)^2_6= t^2 e^{26},\quad
(\Omega^g)^3_4= t^2 (2\, e^{12}- 3\, e^{34}),\quad (\Omega^g)^3_6t^2 e^{36},\quad (\Omega^g)^4_6= t^2 e^{46},
\end{array}
$$
and the non-zero curvature forms $(\Omega^+)^i_j$ of the connection
$\nabla^+$ are
\begin{equation}\label{h3-instanton}
(\Omega^+)^1_2=-(\Omega^+)^3_4=-4t^2(e^{12}-e^{34}).
\end{equation}

Therefore, \eqref{1} is satisfied and the Pontrjagin classes  are
represented by
$$
p_1(\nabla^g)={-3 t^4\over \pi^2} e^{1234}, \quad\quad
p_1(\nabla^+)={-8 t^4\over \pi^2} e^{1234}.
$$

Now, let us consider the new basis $\{f^1,\ldots,f^6\}$ given by
$f^i=e^i$, for $i=1,\ldots,5$, and $f^6=\frac{1}{t} e^6$. In terms
of this basis, the structure equations~\eqref{family-h3} become
$$
d f^1= d f^2= d f^3= d f^4= d f^5 =0, \quad\quad d f^6= -2\, f^{12}
+ 2\, f^{34},
$$
and the family $(J_t,g_t)$ of balanced Hermitian SU($3$)-structures
on $\frh_3$ is given by
\begin{gather*}
J_t f^1=-f^2, \quad J_t f^2=f^1, \quad J_t f^3=-f^4, \quad J_t
f^4=f^3, \quad J_t f^5=-t\, f^6, \quad J_t f^6=\frac{1}{t} f^5,\\
g_t=f^1\otimes f^1+\cdots + f^5\otimes f^5 + t^2 f^6\otimes f^6,
\quad F_t=f^{12}+f^{34}+t\, f^{56}.
\end{gather*}

Let us fix $t'\not= 0$ and denote by $\nabla^+_{t'}$ the connection
corresponding to the balanced structure $(J_{t'},g_{t'})$ in the
previous family. It follows from~\eqref{h3-instanton} that the
non-zero curvature forms $(\Omega^+_{t'})^i_j$ of $\nabla^+_{t'}$
are
$$
(\Omega^+_{t'})^1_2=-(\Omega^+_{t'})^3_4=-4{t'}^2(f^{12}-f^{34}).
$$
Therefore, \eqref{1} and \eqref{2} are satisfied and $\nabla^+_{t'}$
is an SU($3$)-instanton with respect to any other balanced structure
in the family $(J_t,g_t)$.

Let $H(2,1)$ denote the 5-dimensional generalized Heisenberg group,
and let $\Gamma$ be a lattice of maximal rank. The nilpotent Lie
algebra $\frh_3$ is the Lie algebra underlying the compact
nilmanifold $M_3=\Gamma\backslash H(2,1) \times S^1$.

\begin{thrm}\label{solutions-h3-1}
In the notation above, for each $t'\not= t$, we consider the {\rm
SU($3$)}-instanton $\nabla^+_{t'}$. Then we have:
\begin{enumerate}
\item[a)] \hspace{3cm} $dT = {8\pi^2 t^2\over 3t^4 - 8t'^4}
(p_1(\nabla^g_t)-p_1(\nabla^+_{t'}))$,

\medskip

\item[b)] \hspace{3cm} $dT = {\pi^2 t^2\over t^4 - t'^4}
(p_1(\nabla^+_t)-p_1(\nabla^+_{t'}))$.
\end{enumerate}
Hence, for any pair $(t,t')$ such that $8t'^4<3t^4$ we obtain
explicit valid solutions to the heterotic supersymmetry equations
\eqref{sup1} with non-zero flux $H=T$ and constant dilaton
satisfying the three-form Bianchi identity \eqref{acgen} for the
Levi-Civita connection and for the (+)-connection on the compact
nilmanifold $M_3$.

The compact manifold $(M_3,g,J,A=\nabla^+_{t'},R(\nabla^+_t))$
described in b) solves the equations of motion \eqref{mot}.
\end{thrm}

Moreover, we can also use the abelian instanton $A$ given in
Subsection~\ref{card} to find more solutions. In fact, we can take
$dz_1$ and $dz_2$ as (2,0)-forms at the level of the Lie group
$H(2,1) \times \mathbb{R}$ which descend to the forms $\omega^1$ and
$\omega^2$ on the compact nilmanifold $M_3$.

\begin{thrm}\label{solutions-h3-2}
In the notation above and taking $A$ as the abelian {\rm
SU($3$)}-instanton given in \cite{CCDLMZ} we have:
\begin{enumerate}
\item[a)] \hspace{3cm} $dT = {32\pi^2 t^2\over 12 t^4 - 1}
(p_1(\nabla^g_t)-p_1(A))$,

\medskip

\item[b)] \hspace{3cm} $dT = {32\pi^2 t^2\over 32 t^4 - 1}
(p_1(\nabla^+_t)-p_1(A))$.
\end{enumerate}
Thus, for any  $t$ such that $12t^4>1$ we obtain explicit valid
solutions to the heterotic supersymmetry equations \eqref{sup1} with
non-zero flux $H=T$ and constant dilaton satisfying the three-form
Bianchi identity \eqref{acgen} for the Levi-Civita connection and
for the (+)-connection on the compact nilmanifold $M_3$.

The space $(M_3,g,J,A,R(\nabla^+_t))$ described in b)  is a compact
solution to the equations of motion \eqref{mot}.
\end{thrm}

\begin{rmrk}\label{no-solution-h3}
A direct calculation for $\nabla^-$ and for the Chern connection
$\nabla^c$ shows that
$$
p_1(\nabla^-)=0,\quad\quad p_1(\nabla^c)=0.
$$
The nilmanifold $M_3$ is a torus bundle over a complex torus,
therefore we can use the argument given in
Remark~\ref{no-solution-J_0} to conclude that the family above
cannot provide any solution for the connections $\nabla^-$ and
$\nabla^c$.
\end{rmrk}

\section{Balanced Hermitian structures on the Lie algebras $\frh_2$, $\frh_4$ and $\frh_5$}

In this section we construct explicit solutions on compact
nilmanifolds corresponding to the Lie algebras $\frh_2$, $\frh_4$
and $\frh_5$.

Let us consider the complex structure equations
$$
d\omega^1=d\omega^2=0,\quad\quad d\omega^3 \omega^{12}+\omega^{1\bar 1}+b\,\omega^{1\bar2}-\omega^{2\bar2},
$$
where $b\in \mathbb{R}$. Acoording to \cite[Proposition 13]{U}, the
Lie algebras underlying this 1-parameter family of complex equations
are:
\begin{equation}\label{parameter-b}
\mbox{$\frh_2$, for $b\in (-1,1)$;\quad \quad $\frh_4$, for $b=\pm
1$; \quad \quad  $\frh_5$, for any $b$ such that $b^2>1$.}
\end{equation}
Notice that the latter condition defines a 1-parameter family of
complex structures $J$ on the Iwasawa manifold which are not
equivalent to the standard $J_0$.

For any $t\not=0$, let us consider $F$ given by
$$
F=\frac{i}{2} (\omega^{1\bar{1}} + \omega^{2\bar{2}} +
t^2\,\omega^{3\bar{3}}).
$$
Since $D=-1$, $r=s=1$ and the coefficients $u,v,z$ in~\eqref{2forma}
vanish, it follows from Proposition~\ref{balanced}~(ii) that all the
Hermitian structures $(J,F)$ are balanced.

Notice that the associated compact nilmanifolds are $\mathbb{T}^2$
bundles over $\mathbb{T}^4$ for any $b$, whereas the parameter $t$
scales the fiber.

In terms of the real basis of 1-forms $\{e^1,\ldots,e^6\}$ defined
by
$$ e^1+i\, e^2=\omega^1,\quad e^3+i\,
e^4=\omega^2,\quad e^5+i\, e^6=t\, \omega^3,
$$
the structure equations are
\begin{equation}\label{family-h2-h4-h5}
\left\{
  \begin{aligned}
  &d e^1= d e^2= d e^3= d e^4= 0, \\
  &d e^5 = t(b+1)e^{13}+t(b-1)e^{24},\\
  &d e^6= -2t\,e^{12}-t(b-1)e^{14}+t(b+1)e^{23}+2t\, e^{34},
  \end{aligned}
\right.
\end{equation}
the complex structure $J$ is given by $ J e^1=-e^2, J e^3=-e^4, J
e^5=-e^6,$ and the balanced  $J$-Hermitian metric $g=e^1\otimes
e^1+\cdots + e^6\otimes e^6$ has the associated fundamental form $
F=e^{12}+e^{34}+e^{56}. $

Use \eqref{family-h2-h4-h5} to get $dF=2t\, e^{125}+ t(b+1)e^{136}+
t(b-1)e^{145}- t(b+1)e^{235}+ t(b-1)e^{246}- 2t\, e^{345} .$ Due to
\eqref{tor} the torsion $T$ of $\nabla^+$ satisfies

$$\begin{array}{l}
T= -2t\, e^{126}+ t(b-1)e^{135}- t(b+1)e^{146}+ t(b-1)e^{236}+
t(b+1)e^{245}+2t\, e^{346},\\
dT=-4t^2(b^2+3)e^{1234}.\end{array}
$$

A direct calculation using \eqref{pm} gives that the non-zero
curvature forms $(\Omega^+)^i_j$ of the connection $\nabla^+$
 are:

\begin{eqnarray*}
(\Omega^+)^1_2 \!\!\!&=&\!\!\! -4 t^2 e^{12}- 2t^2(b-1)e^{14}+ 2t^2(b+1)e^{23}+ 6t^2 e^{34}-2t^2 b^2e^{56},\\
(\Omega^+)^1_3=(\Omega^+)^2_4 \!\!\!&=&\!\!\! -t^2(b^2+b+1)e^{13}-t^2(b^2-b+1)e^{24},\\
(\Omega^+)^1_4=-(\Omega^+)^2_3 \!\!\!&=&\!\!\! -2t^2 b \, e^{12}-t^2
(b^2-b+1)e^{14}+
t^2(b^2+b+1)e^{23}+2t^2 b\, e^{34}+ 4t^2 b\, e^{56},\\
(\Omega^+)^1_5=(\Omega^+)^2_6 \!\!\!&=&\!\!\! t^2 b\, e^{15}+ t^2 b\, e^{26}-2t^2 e^{46},\\
(\Omega^+)^1_6=-(\Omega^+)^2_5 \!\!\!&=&\!\!\! -t^2 b\, e^{16}+ t^2 b\, e^{25}+2t^2 e^{36},\\
(\Omega^+)^3_4 \!\!\!&=&\!\!\! 6t^2\, e^{12}+2t^2(b-1)e^{14}-2t^2(b+1)e^{23}-4t^2 e^{34}+2t^2 b^2e^{56},\\
(\Omega^+)^3_5=(\Omega^+)^4_6 \!\!\!&=&\!\!\! -2t^2 e^{26}+t^2 b\, e^{35}- t^2 b\, e^{46},\\
(\Omega^+)^3_6=-(\Omega^+)^4_5 \!\!\!&=&\!\!\! 2t^2 e^{16}+t^2 b\, e^{36}+ t^2 b\, e^{45},\\
(\Omega^+)^5_6=-(\Omega^+)^1_2-(\Omega^+)^3_4 \!\!\!&=&\!\!\! -2t^2
e^{12}-2t^2 e^{34}.
\end{eqnarray*}

Similarly, applying \eqref{lc}, we calculate that the non-zero
curvature forms $(\Omega^g)^i_j$ of the Levi-Civita connection
$\nabla^g$ are:

\begin{eqnarray*}
(\Omega^g)^1_2 \!\!\!&=&\!\!\! -3 t^2 e^{12}-\frac32
t^2(b-1)e^{14}+\frac32 t^2(b+1)e^{23}-
\frac{t^2}{2}(b^2-5)e^{34}-\frac{t^2}{2}(b^2+1)e^{56},\\
(\Omega^g)^1_3 \!\!\!&=&\!\!\! -\frac34 t^2(b+1)^2e^{13}-\frac{t^2}{4}(b^2-5)e^{24},\\
(\Omega^g)^1_4 \!\!\!&=&\!\!\! -\frac32 t^2(b-1)e^{12}-\frac34
t^2(b-1)^2e^{14}+\frac{t^2}{4}(b^2-5)e^{23}+
\frac32 t^2(b-1)e^{34}+t^2b\, e^{56},\\
(\Omega^g)^1_5 \!\!\!&=&\!\!\! \frac{t^2}{4}(b+1)^2e^{15}-\frac{t^2}{4}(b-1)^2e^{26}+\frac{t^2}{2}(b-1)e^{46},\\
(\Omega^g)^1_6 \!\!\!&=&\!\!\!
\frac{t^2}{4}(b^2-2b+5)e^{16}+\frac{t^2}{4}(b+1)^2e^{25}+ t^2
e^{36} - \frac{t^2}{2}(b+1)e^{45}, \\
(\Omega^g)^2_3 \!\!\!&=&\!\!\! \frac32 t^2(b+1)e^{12}+\frac{t^2}{4}
(b^2-5)e^{14}-\frac34 t^2
(b+1)^2e^{23}- \frac32 t^2(b+1)e^{34}-t^2b\, e^{56},\\
(\Omega^g)^2_4 \!\!\!&=&\!\!\! -\frac{t^2}{4}(b^2-5)e^{13}-\frac34
t^2(b-1)^2e^{24},\\
(\Omega^g)^2_5 \!\!\!&=&\!\!\!
\frac{t^2}{4}(b+1)^2e^{16}+\frac{t^2}{4}(b-1)^2e^{25}-\frac{t^2}{2}(b+1)e^{36},
\\
(\Omega^g)^2_6 \!\!\!&=&\!\!\!
-\frac{t^2}{4}(b-1)^2e^{15}+\frac{t^2}{4}(b^2+2b+5)e^{26}+\frac{t^2}{2}(b-1)e^{35}-
t^2e^{46},\\
(\Omega^g)^3_4 \!\!\!&=&\!\!\! -\frac{t^2}{2}(b^2-5)e^{12}+\frac32
t^2(b-1)e^{14}-\frac32 t^2(b+1)e^{23}-
3 t^2e^{34}+\frac{t^2}{2}(b^2-1)e^{56},\\
(\Omega^g)^3_5 \!\!\!&=&\!\!\! \frac{t^2}{2}(b-1)e^{26}+\frac{t^2}{4}(b+1)^2e^{35}+\frac{t^2}{4}(b^2-1)e^{46},\\
(\Omega^g)^3_6 \!\!\!&=&\!\!\! t^2
e^{16}-\frac{t^2}{2}(b+1)e^{25}+\frac{t^2}{4}(b^2+2b+5)e^{36}-
\frac{t^2}{4}(b^2-1)e^{45},\\
(\Omega^g)^4_5 \!\!\!&=&\!\!\! -\frac{t^2}{2}(b+1)e^{16}-\frac{t^2}{4}(b^2-1)e^{36}+\frac{t^2}{4}(b-1)^2e^{45},\\
(\Omega^g)^4_6 \!\!\!&=&\!\!\! \frac{t^2}{2}(b-1)e^{15}- t^2
e^{26}+\frac{t^2}{4}(b^2-1)e^{35}+
\frac{t^2}{4}(b^2-2b+5)e^{46},\\
(\Omega^g)^5_6 \!\!\!&=&\!\!\! -\frac{t^2}{2}(b^2+1)e^{12}+ t^2b\,
e^{14}-t^2b\, e^{23}+\frac{t^2}{2}(b^2-1)e^{34}.
\end{eqnarray*}

Hence, $Hol(\nabla^+)\subset {\frak s}{\frak u}(3)$ and the
Pontrjagin classes of the connections $\nabla^g$ and $\nabla^+$ are
represented by
$$
p_1(\nabla^g)= -\frac{t^4}{4\pi^2}(b^4+4b^2+11)e^{1234}, \quad
p_1(\nabla^+)=-\frac{t^4}{\pi^2}(b^4+5b^2+10) e^{1234}.
$$
As we mentioned above, $\frh_5$ is the nilpotent Lie algebra
underlying the Iwasawa manifold. Notice that $\frh_2$ is the Lie
algebra of $H^3\times H^3$, where $H^3$ is the Heisenberg group. Let
us denote by $M_2,M_4,M_5$ any compact nilmanifold whose underlying
Lie algebra is isomorphic to $\frh_2$, $\frh_4$ or $\frh_5$,
respectively. We can take $dz_1$ and $dz_2$ as (2,0)-forms at the
level of the associated Lie group which descend to the forms
$\omega^1$ and $\omega^2$ on $M_2,M_4,M_5$, so using again the
abelian instanton given in Section~\ref{Iwasawa} we get:

\begin{thrm}\label{solutions-h2-h4-h5}
In the notation above and taking $A$ as the abelian {\rm
SU($3$)}-instanton given in \cite{CCDLMZ} we have:
$$
\begin{array}{lll}
& dT =  {16\pi^2 t^2 (b^2+3)\over t^4 (b^4+4b^2+11) - 1}
(p_1(\nabla^g)-p_1(A)),\\[12pt]
& dT =  {16\pi^2 t^2 (b^2+3)\over 4 t^4 (b^4+5b^2+10) - 1}
(p_1(\nabla^+)-p_1(A)).
\end{array}
$$
For any $b\in \mathbb{R}$ we can choose $t\not=0$ such that $$t^4
(b^4+5b^2+10)>1/4 \qquad {\rm and} \qquad  t^4 (b^4+4b^2+11)>1,$$
which, in view of~\eqref{parameter-b}, provides
 explicit valid solutions to the heterotic
supersymmetry equations \eqref{sup1} with non-zero flux $H=T$ and
constant dilaton satisfying the three-form Bianchi identity
\eqref{acgen} for the Levi-Civita connection and for the
(+)-connection on the compact nilmanifolds $M_2,M_4,M_5$.
\end{thrm}

\begin{rmrk}\label{no-solution-family-b}
Finally, a direct calculation for $\nabla^-$ and for the Chern
connection $\nabla^c$ shows that
$$
p_1(\nabla^-)=\frac{t^4}{\pi^2}(b^2+3) e^{1234},\quad\quad
p_1(\nabla^c)=0.
$$
Since $M_2$, $M_4$ and $M_5$ are torus bundles over a complex torus,
notice that the same argument as in Remark~\ref{no-solution-J_0}
shows that the family above cannot provide any satisfactory solution
for the connections $\nabla^-$ and $\nabla^c$.
\end{rmrk}

\section{The space of balanced structures on $\frh_6$}

In this section we study the space of balanced Hermitian structures
on the nilpotent Lie algebra $\frh_6$.

The complex equations
$$
d\omega^1=d\omega^2=0, \quad d\omega^3 = \omega^{12} - \omega^{2\bar
1},
$$
define a complex structure $J$ on $\frh_6$, and any complex
structure on the Lie algebra $\frh_6$ is equivalent to
$J$~\cite[Corollary 15]{U}. Moreover, it is easy to see that any
$J$-balanced structure $F$ is equivalent to one of the form
$$
F=\frac{i}{2} (\omega^{1\bar{1}} + \omega^{2\bar{2}} +
t^2\,\omega^{3\bar{3}}),
$$
for some $t\not=0$.

From a real point of view, the whole space of balanced Hermitian
structures on $\frh_6$ is described as follows. Let us consider the
basis of 1-forms $\{e^1,\ldots,e^6\}$ given by
$$
e^1+i\, e^2=\omega^1,\quad e^3+i\, e^4=\omega^2,\quad e^5+i\,
e^6=t\, \omega^3.
$$
Now, in terms of this basis, we have the structure equations
\begin{equation}\label{family-h6}
\left\{
  \begin{aligned}
  &d e^1= d e^2= d e^3= d e^4 =0, \\
  &d e^5= 2t\, e^{13},\\
  &d e^6= 2t\, e^{14}.
  \end{aligned}
\right.
\end{equation}
The complex structure $J$ is given by $ J e^1=-e^2, J e^3=-e^4, d J
e^5=-e^6,$ the $J$-Hermitian metric $g=e^1\otimes e^1+\cdots +
e^6\otimes e^6$ has the associated fundamental form
$F=e^{12}+e^{34}+e^{56}.$

The structure equations \eqref{family-h6} yield $dF=2t (e^{136} -
e^{145})$. Consequently, applying \eqref{tor}, we obtain that the
torsion $T$ of $\nabla^+$ satisfies
$$
T= -2t(e^{236}-e^{245}), \qquad dT=-8t^2 e^{1234}.
$$

Using \eqref{pm} we calculate that the non-zero curvature forms
$(\Omega^+)^i_j$ for the connection $\nabla^+$ are given by:

$$
\begin{array}{l}
(\Omega^+)^1_2= 2 t^2(e^{34}+e^{56}),\quad
(\Omega^+)^1_3=(\Omega^+)^2_4=-t^2(3e^{13}+e^{24}), \quad
(\Omega^+)^1_4=-(\Omega^+)^2_3=-t^2(3e^{14}-e^{23}), \\[10pt]
(\Omega^+)^1_5=(\Omega^+)^2_6= t^2(e^{15}-e^{26}),\quad
(\Omega^+)^1_6=-(\Omega^+)^2_5= t^2(e^{16}+e^{25}), \quad
(\Omega^+)^3_4= 2t^2(e^{12}-e^{56}), \\[10pt]
(\Omega^+)^3_5=(\Omega^+)^4_6= t^2(e^{35}+e^{46}), \quad
(\Omega^+)^3_6=-(\Omega^+)^4_5= -t^2(e^{36}-e^{45}), \quad
(\Omega^+)^5_6= -2t^2(e^{12}+e^{34}),
\end{array}
$$
so \eqref{1} holds and the first Pontrjagin class is represented by
$$
p_1(\nabla^+)=-{2t^4\over \pi^2} e^{1234}.
$$
Let us denote by $M_6$ any compact nilmanifold whose underlying Lie
algebra is isomorphic to $\frh_6$. We can take $dz_1$ and $dz_2$ as
(2,0)-forms at the level of the Lie group corresponding to $\frh_6$
which descend to the forms $\omega^1$ and $\omega^2$ on $M_6$, so
using again the abelian instanton given in Section~\ref{Iwasawa} we
get:

\begin{thrm}\label{solutions-h6}
In the notation above and taking $A$ as the abelian {\rm
SU($3$)}-instanton given in \cite{CCDLMZ} we have:
$$
\begin{array}{lll}
& dT =  {32\pi^2 t^2 \over 8 t^4 - 1} (p_1(\nabla^+)-p_1(A)).
\end{array}
$$
Thus, for any  $t$ such that $t^4>\frac{1}{8}$ we obtain explicit
valid solutions to the heterotic supersymmetry equations
\eqref{sup1} with non-zero flux $H=T$ and constant dilaton
satisfying the three-form Bianchi identity \eqref{acgen} for  the
(+)-connection on the compact nilmanifold $M_6$.
\end{thrm}

\begin{rmrk}\label{no-solution-h6}
The Pontrjagin classes of the Levi-Civita connection, $\nabla^-$ and
the Chern connection are represented by
$$
p_1(\nabla^g)= 0,\quad\quad p_1(\nabla^-)= {2t^4\over \pi^2}
e^{1234},\quad\quad p_1(\nabla^c)=0.
$$
Since the nilmanifold $M_6$ is a torus bundle over a complex torus,
the same argument as in Remark~\ref{no-solution-J_0} shows that
there is no way to find a satisfactory solution for the connections
$\nabla^g$, $\nabla^-$ and $\nabla^c$ on the whole space of
invariant balanced Hermitian structures on $M$.
\end{rmrk}

\section{Balanced structures on the Lie algebra $\frh^-_{19}$}
In this section we construct compact valid solutions to \eqref{sup1}
with non-zero flux and constant dilaton satisfying anomaly
cancellation condition \eqref{acgen} using the curvature $R^c$ of
the Chern connection.

Consider the complex structure equations
$$
d\omega^1=0,\quad d\omega^2=\omega^{13}+\omega^{1\bar3},\quad
d\omega^3= i(\omega^{1\bar2}-\omega^{2\bar1}),
$$
which in view of \eqref{nonilp-struct} correspond to a complex
structure $J$ on the $3$-step nilpotent Lie algebra $\frh^-_{19}$.

The associated real structure equations are
\begin{equation}\label{equations-h19}
\left\{
  \begin{aligned}
  &de^1=de^2=de^5=0, \\
  &de^3=2e^{15},\\
  &de^{4}=2e^{25},\\
  &de^6=2(e^{13}+e^{24}),
  \end{aligned}
\right.
\end{equation}
and the complex structure $J$ is given by $Je^1=-e^2, Je^3=-e^4,
Je^5=-e^6.$ The fundamental form $F$ of the $J$-Hermitian metric
$g=e^1\otimes e^1+\cdots+e^6\otimes e^6$ is given by
$F=e^{12}+e^{34}+e^{56}.$ It follows from
Proposition~\ref{balanced}~(i) that the structure $(J,g)$ is
balanced.

The structure equations \eqref{equations-h19} imply
$dF=-2(e^{135}+e^{145}-e^{235}+e^{245})$. Apply \eqref{tor} to
verify that the torsion $T$ satisfies

$$T=2(e^{136}+e^{146}-e^{236}+e^{246})\qquad dT=-8(e^{1234}+e^{1256}).
$$

Using \eqref{lc}, \eqref{pm} and \eqref{c} we obtain that the
non-zero curvature forms  $(\Omega^c)^i_j$ and $(\Omega^+)^i_j$ of
the Chern connection and  the (+)-connection  are given by:

$$
\begin{array}{l}
(\Omega^c)^1_2=-2e^{34}-2e^{56},\quad
(\Omega^c)^1_3=(\Omega^c)^2_4=-e^{13}-e^{24},\quad
(\Omega^c)^1_4=-(\Omega^c)^2_3=2e^{13}+e^{14}-e^{23}+2e^{24},\\[5pt]
(\Omega^c)^1_5=(\Omega^c)^2_6=e^{16}-e^{25},\quad(\Omega^c)^1_6=(\Omega^c)^2_5-e^{15}-e^{26},\quad (\Omega^c)^3_4=-2e^{12}+2e^{56},\\[5pt]
(\Omega^c)^3_5=(\Omega^c)^4_6=-e^{36}+e^{45},\quad
(\Omega^c)^3_6=-(\Omega^c)^4_5=e^{35}+e^{46},\quad
(\Omega^c)^5_6=-(\Omega^c)^1_2-(\Omega^c)^3_4=2e^{12}+2e^{34}.
\end{array}
$$

$$
\begin{array}{l}
(\Omega^+)^1_2=-2e^{34}+2e^{56},\quad(\Omega^+)^1_3=(\Omega^+)^2_4=-3e^{13}-3e^{24},\\[5pt]
(\Omega^+)^1_4=-(\Omega^+)^2_3=-2e^{13}-e^{14}+e^{23}-2e^{24},\quad
(\Omega^+)^1_5=(\Omega^+)^2_6=-3e^{15}-2e^{16}-e^{26},\\[5pt]
(\Omega^+)^1_6=-(\Omega^+)^2_5=-e^{16}+3e^{25}+2e^{26},\quad
(\Omega^+)^3_4=-2e^{12}-2e^{56},\\[5pt]
(\Omega^+)^3_5=(\Omega^+)^4_6=e^{35}+2e^{36}-e^{46},\quad(\Omega^+)^3_6=-(\Omega^+)^4_5=-e^{36}-e^{45}-2e^{46},\\[5pt]
(\Omega^+)^5_6=-(\Omega^+)^1_2-(\Omega^+)^3_4=2e^{12}+2e^{34};
\end{array}
$$

A direct calculation shows that  the Pontrjagin classes  are
represented by
$$ p_1(\nabla^+)=-\frac{2}{\pi^2}(3e^{1234}+e^{1256}),
\quad\quad p_1(\nabla^c)=-\frac{2}{\pi^2}(e^{1234}+e^{1256}).
$$

Let $M_{19}$ be a compact nilmanifold corresponding to the Lie
algebra $\frh^-_{19}$. From \eqref{equations-h19} we have that
$M_{19}$ is an $S^1$-bundle over a compact 5-nilmanifold $N$, which
is a $\mathbb{T}^2$-bundle over $\mathbb{T}^3$.

\begin{lemma}\label{instanton-h19}
For each $\lambda,\mu \in \mathbb{R}$, let $A_{\lambda,\mu}$ be the
{\rm U($3$)}-connection on $M_{19}$ with respect to structure
$(J,g)$ defined by the connection forms
$$
(\sigma^{A_{\lambda,\mu}})^2_3=(\sigma^{A_{\lambda,\mu}})^2_5=(\sigma^{A_{\lambda,\mu}})^4_5-\lambda\, e^1 - \mu\, e^6,\quad\quad
(\sigma^{A_{\lambda,\mu}})^i_j= \lambda\, e^1 + \mu\, e^6,
$$
for $1\leq i< j\leq 6$ such that $(i,j)\not=(2,3),(2,5),(4,5)$.
Then, $A_{\lambda,\mu}$ is an {\rm SU($3$)}-instanton and
$$p_1(A_{\lambda,\mu})=-\frac{15}{\pi^2} \mu^2 e^{1234}.$$
\end{lemma}

\begin{proof}
A direct calculation shows that the curvature forms
$(\Omega^{A_{\lambda,\mu}})^i_j$ of the connection $A_{\lambda,\mu}$
are given by
$$(\Omega^{A_{\lambda,\mu}})^2_3=(\Omega^{A_{\lambda,\mu}})^2_5= (\Omega^{A_{\lambda,\mu}})^4_5= -2\mu(e^{13}+e^{24}),\quad\quad
(\Omega^{A_{\lambda,\mu}})^i_j= 2\mu(e^{13}+e^{24}),$$ for $1\leq i<
j\leq 6$ such that $(i,j)\not=(2,3),(2,5),(4,5)$. Now it is clear
that $A_{\lambda,\mu}$ satisfies \eqref{2}.
\end{proof}

\begin{thrm}\label{solution-h19}
Let $A_{\lambda,\mu}$ be the {\rm SU($3$)}-instanton above.
\begin{enumerate}
\item[{\rm (i)}] If $\mu^2= \frac{4}{15}$, then
$$
dT = 4\pi^2 (p_1(\nabla^+)-p_1(A_{\lambda,\mu})).$$ \item[{\rm
(ii)}] If $\mu=0$, then $p_1(A_{\lambda,0})=0$ and
$$
dT = 4\pi^2\, (p_1(\nabla^c)-p_1(A_{\lambda,0})).
$$
\end{enumerate}
Hence, we obtain explicit valid solutions to the heterotic
supersymmetry equations \eqref{sup1} with non-zero flux $H=T$ and
constant dilaton satisfying the three-form Bianchi identity
\eqref{acgen} for the Chern connection and the (+)-connection on the
compact nilmanifold $M_{19}$.
\end{thrm}

\begin{rmrk}
During the preparation of the paper we learned that a compact
example solving \eqref{sup1} with non-zero flux, constant dilaton
satisfying \eqref{acgen} with respect to a metric connection on the
tangent bundle, and trivial instanton ($A=0$) on $M_3$ is announced
\cite{GrP}.
\end{rmrk}

\medskip
\noindent {\bf Acknowledgments.}  We thank George Papadopoulos for
very useful discussions. This work has been partially supported
through grant MEC (Spain) MTM2005-08757-C04-02. S.I. is partially
supported by the Contract 154/2008 with the University of Sofia
`St.Kl.Ohridski`. S.I. is a Senior Associate to the Abdus Salam ICTP, Trieste.


\begin{thebibliography}{33}

\bibitem{AI} B. Alexandrov, S. Ivanov, {\em Vanishing theorems on Hermitian manifolds},
 Diff. Geom. Appl. {\bf 14} (3) (2001), 251--265.

\bibitem {BBDG} K. Becker, M. Becker, K. Dasgupta, P.S. Green,
{\em Compactifications of Heterotic Theory on Non-K\"ahler Complex
Manifolds: I}, JHEP 0304 (2003) 007.

\bibitem{BBE}  K. Becker, M. Becker, K. Dasgupta, P.S. Green, E. Sharpe,
{\em Compactifications of Heterotic Strings on Non-K\"ahler Complex
Manifolds: II}, Nucl. Phys. {\bf B} {\bf 678} (2004), 19--100.

\bibitem{BBDP} K. Becker, M. Becker, K. Dasgupta, S. Prokushkin,
{\em Properties from heterotic vacua from superpotentials},
hep-th/0304001.

\bibitem{y4} K. Becker, M. Becker, J-X. Fu, L-S. Tseng, S-T. Yau,
{\em Anomaly Cancellation and Smooth Non-Kahler Solutions in
Heterotic String Theory}, Nucl. Phys. {\bf B} {\bf 751} (2006),
108--128.

\bibitem{Berg} E.A. Bergshoeff, M. de Roo, {\em The quartic
effective action of the heterotic string and supersymmetry}, Nucl.
Phys. {\bf B} {\bf 328} (1989), 439.

\bibitem{Bis} J.-M. Bismut, {\em A local index theorem for non-K\"ahler
manifolds}, Math. Ann. {\bf 284} (1989), 681--699.

\bibitem{Car1}  G.L. Cardoso, G. Curio, G. Dall'Agata, D. Lust, {\em
BPS Action and Superpotential for Heterotic String Compactifications
with Fluxes}, JHEP 0310 (2003) 004.

\bibitem{CCDLMZ} G.L. Cardoso, G. Curio, G. Dall'Agata, D. Lust, P. Manousselis, G. Zoupanos,
{\em Non-K\"aeler string back-grounds and their five torsion
classes}, Nucl. Phys. {\bf B} {\bf 652} (2003), 5--34.

\bibitem {CS} S. Chiossi, S. Salamon, {\em The intrinsic torsion of SU(3) and
  $G_2$-structures}, Differential Geometry, Valencia 2001, World Sci.
  Publishing, 2002, pp. 115--133.

\bibitem{CFGU} L.A. Cordero, M. Fern\'andez, A. Gray, L. Ugarte, {\em Compact nilmanifolds
with nilpotent complex structure: Dolbeault cohomology}, Trans.
Amer. Math. Soc.  {\bf 352} (2000),  5405--5433.

\bibitem{CDev} E. Corrigan, C. Devchand, D.B. Fairlie, J. Nuyts, {\em First-order
equations for gauge fields in spaces of dimension greater than
four}, Nucl. Phys. {\bf B} {\bf 214} (1983),  no. 3, 452--464.

\bibitem{DRS} K. Dasgupta, G. Rajesh, S. Sethi, {\em M theory,
orientifolds and G-flux}, JHEP {\bf 0211}, 006 (2002).

\bibitem{DFG} K. Dasgupta, H. Firouzjahi, R. Gwyn, {\em On the warped heterotic axion}, 
arXiv:0803.3828 [hep-th], to appear in JHEP.

\bibitem{Bwit} B. de Wit, D.J. Smit, N.D. Hari Dass, {\em Residual Supersimmetry
Of Compactified D=10 Supergravity}, Nucl. Phys. {\bf B} {\bf 283}
(1987), 165.

\bibitem{FPS} A. Fino, M. Parton, S. Salamon, {\em Families of strong
KT structures in six dimensions}, Comment. Math. Helv.  {\bf 79}
(2004), 317--340.

\bibitem{FGW} D.Z. Freedman, G.W. Gibbons, P.C. West, {\em Ten Into Four Won't Go},
Phys. Lett. {\bf B 124} (1983), 491.

\bibitem {FI} Th. Friedrich, S. Ivanov, {\em Parallel spinors and connections
  with skew-symmetric torsion in string theory}, Asian J. Math. {\bf 6} (2002), 303--335.

\bibitem{y2} J-X. Fu, S-T. Yau, {\em Existence of Supersymmetric Hermitian Metrics
with Torsion on Non-Kaehler Manifolds}, arXiv:hep-th/0509028.

\bibitem{y3} J-X. Fu, S-T. Yau, {\em The theory of superstring with flux on non-K\"ahler
manifolds and the complex Monge-Ampere equation},
arXiv:hep-th/0604063.

\bibitem {GKMW} J. Gauntlett, N. Kim, D. Martelli, D. Waldram, {\em Fivebranes
  wrapped on SLAG three-cycles and related geometry}, JHEP 0111 (2001) 018.

\bibitem {GMPW} J.P. Gauntlett, D. Martelli, S. Pakis, D. Waldram,
{\em  G-Structures and Wrapped NS5-Branes}, Commun. Math. Phys. {\bf
247} (2004), 421--445.

\bibitem {GMW} J. Gauntlett, D. Martelli, D. Waldram, {\em Superstrings with
  Intrinsic torsion}, Phys. Rev. {\bf D69} (2004) 086002.

\bibitem{GPap} J. Gillard, G. Papadopoulos, D. Tsimpis, {\em
Anomaly, Fluxes and (2,0) Heterotic-String Compactifications}, JHEP
0306 (2003) 035.

\bibitem{GP} E. Goldstein, S. Prokushkin, {\em Geometric Model for Complex
Non-K\"aehler Manifolds with SU(3) Structure}, Commun. Math. Phys.
{\bf 251} (2004), 65--78.

\bibitem{GrP} G. Grantcharov, Y-S. Poon, Talk of G. Grantcharov in
 the Workshop 'Special Geometries in Mathematical Physics',
 Kulungsborn, March 30-April 4, 2008.

\bibitem{GrH} A. Gray, L. Hervella, {\em The sixteen classes of
almost Hermitian manifolds and their linear invariants}, Ann. Mat.
Pura Appl. (4) {\bf 123} (1980), 35--58.

\bibitem{HP1} P.S. Howe, G. Papadopoulos, {\em Ultraviolet
behavior of two-dimensional supersymmetric non-linear sigma models},
Nucl. Phys. {\bf B} {\bf 289} (1987), 264.

\bibitem{Hull} C.M. Hull, {\em Anomalies, ambiquities and
superstrings}, Phys. Lett. {\bf B 167} (1986), 51.

\bibitem{HT} C.M. Hull, P.K. Townsend, {\em The two loop beta
function for sigma models with torsion}, Phys. Lett. {\bf B 191}
(1987), 115.

\bibitem{HuW} C.M. Hull, E. Witten, {\em Supersymmetric sigma
models and the Heterotic String}, Phys. Lett. {\bf B 160} (1985),
398.

\bibitem{II} P. Ivanov, S. Ivanov, {\em SU(3)-instantons and
$G_2,Spin(7)$-Heterotic string solitons}, Comm. Math. Phys. {\bf
259} (2005), 79--102.

\bibitem {IP2} S. Ivanov, G. Papadopoulos, {\em A no-go theorem for string warped compactifications},
 Phys. Lett. {\bf B 497} (2001), 309--316.

\bibitem {IP1} S. Ivanov, G. Papadopoulos, {\em Vanishing Theorems and String Backgrounds},
 Class. Quant. Grav. {\bf 18} (2001), 1089--1110.
 
 \bibitem{KY} T. Kimura, P. Yi, {\em Comments on heterotic flux compactifications}, JHEP {\bf 0607}, 030 
 (2006) [arXiv:hep-th/0605247].

\bibitem{y1} J. Li, S-T. Yau, {\em The Existence of Supersymmetric String Theory with Torsion},
 J. Diff. Geom. {\bf 70} (2005), no. 1, 143--181.

\bibitem{Mi} M.L. Michelsohn, {\em On the existence of special metrics
in complex geometry}, Acta Math. {\bf 149} (1982), no. 3-4,
261--295.

\bibitem{Sen} A. Sen, {\em (2,0) supersymmetry and space-time
supersymmetry in the heterotic strin theory}, Nucl. Phys. {\bf B}
{\bf 167} (1986), 289.

\bibitem {Str} A. Strominger, {\em Superstrings with torsion}, Nucl. Phys.
{\bf B} {\bf 274} (1986), 253.

\bibitem{U} L. Ugarte,
{\em Hermitian structures on six-dimensional nilmanifolds},
Transform. Groups {\bf 12} (2007), 175--202.


\end{thebibliography}
\end{document}